\documentclass[12pt]{article}
\usepackage{epsf}
\usepackage[usenames]{color}
\usepackage{epsfig,amsmath,times,amssymb,euscript,graphicx,cite,amssymb}
\usepackage{pstricks,fancybox,hyperref}
\fontfamily{pnc}\selectfont

\definecolor{myorange}{rgb}{1,0.5,0}
\begin{document}
\setlength{\parskip}{2ex} \setlength{\parindent}{0em}
\setlength{\baselineskip}{3ex}
\newcommand{\onefigure}[2]{\begin{figure}[htbp]
         \caption{\small #2\label{#1}(#1)}
         \end{figure}}
\newcommand{\onefigurenocap}[1]{\begin{figure}[h]
         \begin{center}\leavevmode\epsfbox{#1.eps}\end{center}
         \end{figure}}
\renewcommand{\onefigure}[2]{\begin{figure}[htbp]
         \begin{center}\leavevmode\epsfbox{#1.eps}\end{center}
         \caption{\small #2\label{#1}}
         \end{figure}}
\newcommand{\comment}[1]{}
\newcommand{\myref}[1]{(\ref{#1})}
\newcommand{\secref}[1]{sec.~\protect\ref{#1}}
\newcommand{\figref}[1]{Fig.~\protect\ref{#1}}
\def\sl2z{SL(2,\Z)}
\newcommand{\be}{\begin{equation}}
\newcommand{\ee}{\end{equation}}
\newcommand{\bea}{\begin{eqnarray}}
\newcommand{\eea}{\end{eqnarray}}
\newcommand{\nn}{\nonumber}
\newcommand{\unit}{1\!\!1}
\newcommand{\mt}{\widehat{t}}
\newcommand{\R}{\bf R}
\newcommand{\X}{{\bf X}}
\newcommand{\T}{{\bf T}}
\newcommand{\PP}{\bf P}
\newcommand{\CC}{\bf C}
\newcommand {\us}[1]{\underline{s_{#1}}}
\newcommand {\uk}{\underline{k}}
\newcommand{\me}{\mathellipsis}
\newdimen\tableauside\tableauside=1.0ex
\newdimen\tableaurule\tableaurule=0.4pt
\newdimen\tableaustep
\def\phantomhrule#1{\hbox{\vbox to0pt{\hrule height\tableaurule width#1\vss}}}
\def\phantomvrule#1{\vbox{\hbox to0pt{\vrule width\tableaurule height#1\hss}}}
\def\sqr{\vbox{%
\phantomhrule\tableaustep
\hbox{\phantomvrule\tableaustep\kern\tableaustep\phantomvrule\tableaustep}%
\hbox{\vbox{\phantomhrule\tableauside}\kern-\tableaurule}}}
\def\squares#1{\hbox{\count0=#1\noindent\loop\sqr
\advance\count0 by-1 \ifnum\count0>0\repeat}}
\def\tableau#1{\vcenter{\offinterlineskip
\tableaustep=\tableauside\advance\tableaustep by-\tableaurule
\kern\normallineskip\hbox
    {\kern\normallineskip\vbox
      {\gettableau#1 0 }%
     \kern\normallineskip\kern\tableaurule}%
  \kern\normallineskip\kern\tableaurule}}
\def\gettableau#1 {\ifnum#1=0\let\next=\null\else
  \squares{#1}\let\next=\gettableau\fi\next}

\tableauside=1.0ex \tableaurule=0.4pt

\bibliographystyle{utphys}
\setcounter{page}{1} \pagestyle{plain}
\numberwithin{equation}{section}

\begin{titlepage}
\begin{center}
 \hfill\\ \vskip 1cm { \large Generalizations of  Nekrasov-Okounkov Identity} \vskip 0.5cm
{ Amer\,\,Iqbal$^{1,2}$ \,\,\,Shaheen Nazir$^{3}$\,\,\, Zahid Raza$^{3}$\,\,\,Zain Saleem$^{4}$}\\
\vskip 0.5cm
$^{1}${ Department of Physics,\\
LUMS School of Science and Engineering,\\
U-Block, D.H.A, Lahore, Pakistan.\\} \vskip 0.5cm
$^{2}${Abdus Salam School of Mathematical Sciences,\\
G.C. University, Lahore, Pakistan.\\} \vskip 0.5cm
$^{3}${Department of Mathematics,\\
National University of Computer and Emerging Sciences,\\ B-Block, Faisal Town, Lahore, Pakistan.\\} \vskip 0.5cm
$^{4}$ {National Center for Physics\\
Quaid-e-Azam University, Islamabad, Pakistan.}
\end{center}
\vskip 2 cm
\begin{abstract}
Nekrasov-Okounkov identity gives a product representation of the sum over partitions of a certain function of partition hook length. In this paper we give several generalizations of the Nekrasov-Okounkov identity using the cyclic symmetry of the topological vertex.
\end{abstract}
\end{titlepage}

\section{Introduction}
The Nekrasov-Okounkov identity was discovered by Nekrasov and Okounkov while studying the supersymmetric gauge theories \cite{NO}. It relates the powers of Euler products to the sum over product of partition hook lengths:
\newline
\begin{equation}
\label{NO}
{\sum_{\lambda}z^{|\lambda|}\,\,\prod_{(i,j)\in \lambda}\frac{%
h(i,j)^2-t^2}{h(i,j)^2}}={\prod_{k\geq 1}(1-z^k)^{t^2-1}}
\end{equation}%
\newline
where ${\lambda} $ is the partition of $n\in \mathbb{Z}^{+}\bigcup\{0\}$ and
we denote by $\lambda_{i}$ the parts of $\lambda$ such that $\lambda_{1}\geq
\lambda_{2}\geq \lambda_{3}\geq \cdots$. $\ell(\lambda)$ denotes the number
of non-zero parts of $\lambda$ and $|\lambda|=\sum_{i}\lambda_{i}$ denotes
the size of the partition. We can represent this partition as a Young
diagram such that there are $\lambda_{i}$ boxes in the $i$-th column. We denote by $\lambda$ both
the partition and its Young diagram. The box in a Young diagram has
coordinates $(i,j)\in \lambda$ if $1\leq i\leq \ell(\lambda),\,1\leq j\leq
\lambda_{i}$. We denote the transpose of the partition $\lambda$ by $%
\lambda^{t}$ which is obtained by interchanging the columns and the rows of $%
\lambda$. The arm length $a(i,j)$, the leg length $\ell(i,j)$ and the hook
length $h(i,j)$ are given by
\begin{eqnarray}
a_{\lambda}(i,j)&=&\lambda^{t}_{j}-i\,,\,\,\, \ell_{\lambda}(i,j)=\lambda_{i}-j \\\nn
h_{\lambda}(i,j)&=&a(i,j)+\ell(i,j)+1=\lambda_{i}+\lambda^{t}_{j}-i-j+1
\end{eqnarray}

\begin{figure}[h]\begin{center}
$\begin{array}{c@{\hspace{1in}}c} \multicolumn{1}{l}{\mbox{}} &
    \multicolumn{1}{l}{\mbox{}} \\ [-0.53cm]
{\begin{pspicture}(7,0)(5,5)

\psframe[unit=0.5cm, linecolor=white, fillstyle=solid,
fillcolor=lightgray](3,3)(4,6)

\psframe[unit=0.5cm, linecolor=white, fillstyle=solid,
fillcolor=lightgray](4,3)(6,5)

\psframe[unit=0.5cm, linecolor=white, fillstyle=solid,
fillcolor=lightgray](6,3)(7,4)

\psframe[unit=0.5cm, linecolor=white, fillstyle=solid,
fillcolor=lightgray](2,3)(3,8)


 \psgrid[unit=0.5cm, subgriddiv=1,
gridcolor=myorange, %
gridlabelcolor=white]%
(2,3)(9,10)

\psframe[unit=0.5cm, linecolor=white, fillstyle=solid,
fillcolor=lightgray](17,3)(19,6)

\psframe[unit=0.5cm, linecolor=white, fillstyle=solid,
fillcolor=lightgray](17,7)(18,8)

\psframe[unit=0.5cm, linecolor=white, fillstyle=solid,
fillcolor=lightgray](19,3)(21,5)

\psframe[unit=0.5cm, linecolor=white, fillstyle=solid,
fillcolor=lightgray](21,3)(22,4)

\psframe[unit=0.5cm, linecolor=white, fillstyle=solid,
fillcolor=lightgray](17,5)(18,7)

\psframe[unit=0.5cm, linecolor=white,
fillstyle=crosshatch,hatchwidth=0.3pt,hatchsep=2.5pt,
fillcolor=lightgray](18,4)(19,5)

\psline[unit=0.5cm, linecolor=red,linewidth=2pt](18.5,5)(18.5,6)
\psline[unit=0.5cm, linecolor=red,linewidth=2pt,linestyle=dashed](19,4.5)(21,4.5)

\psgrid[unit=0.5cm, subgriddiv=1,
gridcolor=myorange, %
gridlabelcolor=white]%
(17,3)(24,10)

\put(1,1){$\lambda_{1}$}
\put(1.55,1){$\lambda_{2}$}
\put(2.05,1){$\lambda_{3}$}
\put(2.5,1){$\lambda_{4}$}
\put(3,1){$\lambda_{5}$}
\put(4,1){$\rightarrow\,i$}
\put(0.45,1.55){$\lambda^{t}_{1}$}
\put(0.45,2.1){$\lambda^{t}_{2}$}
\put(0.45,2.6){$\lambda^{t}_{3}$}
\put(0.45,3.1){$\lambda^{t}_{4}$}
\put(0.45,3.6){$\lambda^{t}_{5}$}
\put(0.45,4.5){$\uparrow $}
\put(0.45,4.9){$j $}
\put(2.5,0.5){$(a)$}
\put(10,0.5){$(b)$}

\end{pspicture}}
\\ [-1.2cm] \mbox{} & \mbox{}
\end{array}$
\caption{(a) Young diagram for $\lambda=\{5,3,2,2,1\}$. The transpose of this is $\lambda^{t}=\{5,4,2,1,1\}$. (b) The horizontal dotted line is the arm length and the vertical line is the leg length of the box with coordinates $(2,2)$.}
\label{f7}\end{center}
\end{figure}
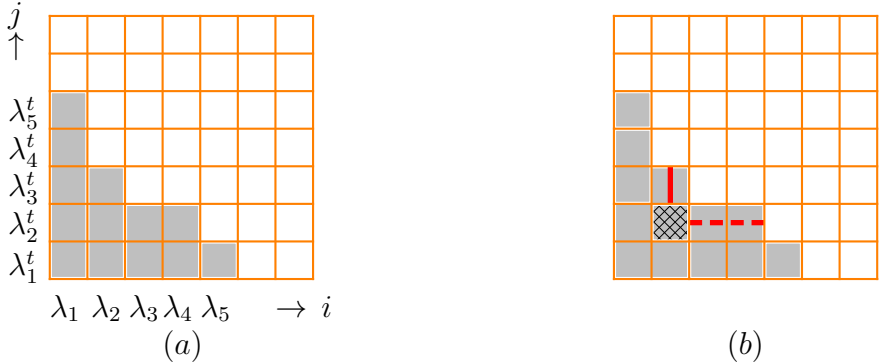

In an earlier paper \cite{adjoint} we stated a two-parameter generalization of NO-identity Eq(\ref{NO}):
\bea
\sum_{\lambda}x^{|\lambda|}\,\,\prod_{s\in \lambda}\frac{
\Big(a(s)+1+\vartheta\,\ell(s)-t\Big)\Big(a(s)+\vartheta(\ell(s)+1)+t\Big)}{\Big(a(s)+1+\vartheta\,\ell(s)\Big)
\Big(a(s)+\vartheta(\ell(s)+1)\Big)}=
\prod_{k\geq 1}(1-x^k)^{(t-1)(\frac{t}{\vartheta}+1)}
\eea
This identity also comes from the study of the Seiberg-Witten theory and its relation with random partitions \cite{adjoint}. In this paper we will give other generalizations of Nekrasov-Okounkov identity based on the topological vertex formalism \cite{mmk,tv,Okounkov:2003sp} for computing the partition function of the Seiberg-Witten theory. We will not discuss the physics behind these identities and will therefore not discuss the Seiberg-Witten theory but will rather use the topological vertex and its symmetries as a tool to generate such identities.

The topological vertex is a function of a complex parameter $q$, indexed by three partitions and defined in terms of the Schur and skew-Schur functions, see \cite{Okounkov:2003sp}\footnote{We will use the definition of the topological vertex given in \cite{Okounkov:2003sp} which differs slightly from the one given in \cite{tv}.} :
\begin{equation}
C_{\lambda\, \mu\, \nu }(q)=q^{\frac{{k(\mu )}}{2}}s_{\nu ^{t}}(q^{-\rho
})\sum_{\eta }s_{\lambda ^{t}/\eta }(q^{-\rho -\nu })s_{\mu /\eta }(q^{-\rho
-\nu ^{t}})\,
\end{equation}
where $s_{\nu/\eta}(x_{1},x_{2},\cdots)$ is the skew-Schur function, $q^{-\rho-\nu}=\{q^{\frac{1}{2}-\nu_{1}},q^{\frac{3}{2}-\nu_{2}},q^{\frac{5}{2}-\nu_{3}},\cdots\}$ and $\kappa(\nu)=\sum_{i}\nu_{i}(\nu_{i}+1-i)$. Notice that $\lambda$ and $\mu$ appear to be almost on the same footing but $\nu$ is treated very differently. An important property of the topological vertex, which is not obvious from the definition given above, is that it is cyclically symmetric:
\bea
C_{\lambda\,\mu\,\nu}(q)=C_{\mu\,\nu\,\lambda}(q)=C_{\nu\,\lambda\,\mu}(q)
\eea
This cyclic symmetry implies interesting identities involving Schur functions,
\bea
s_{\mu}(q^{-\rho})s_{\lambda}(q^{-\rho-\mu})=q^{-\frac{\kappa(\mu)+\kappa(\lambda)}{2}}\,\sum_{\eta}s_{\lambda^t/\eta}(q^{-\rho})\,s_{\mu^t/\eta}(q^{-\rho})
\eea
The topological vertex is a combinatorial object and can be interpreted as the generating function of plane partitions (3D partitions) with certain constraints. This combinatorial definition makes the cyclic symmetry manifest as discussed in \cite{Okounkov:2003sp}.

\section{Derivation of Nekrasov-Okounkov Identity}
In this section we will derive the Nekrasov-Okounkov identity. The method of generalization this identity will also discussed.

Consider the following oriented trivalent graph:

\begin{pspicture}(-2,0)(3,3)
\label{graph1}
\psline[linecolor=black,linewidth=2pt]{->}(5,1.5)(6,2.5)
\psline[linecolor=black,linewidth=2pt]{->}(6,2.5)(6,3.2)
\psline[linecolor=black,linewidth=2pt]{->}(6,2.5)(6.7,2.5)
\psline[linecolor=black,linewidth=2pt]{->}(4.3,1.5)(5,1.5)
\psline[linecolor=black,linewidth=2pt]{->}(5,0.8)(5,1.5)
\put(6.5,2.6){$\emptyset$}
\put(5.7,2.9){$\emptyset$}
\put(5.6,1.7){$\lambda$}
\put(5.1,0.8){$\emptyset$}
\put(4.3,1.1){$\emptyset$}

\end{pspicture}

We will associate with each edge of the graph a partition $\lambda$ and a factor $z^{|\lambda|}$ such that the empty partition $\emptyset$ is associated with the external edges. With each vertex we will associate the ordered triplet $\{\lambda, \mu,\nu\}$ and a factor $C_{\lambda\,\mu\,\nu}(q)$ where $\lambda,\mu$ and $\nu$ are the partitions associated with the edges of that vertex such that the edges are directed into the vertex. If an edge is outgoing from a vertex we replace the partition $\lambda$ associated with that edge with $\lambda^{t}$. We take the convention that the order is given by anticlockwise going around the vertex. We will see later that since the topological vertex is cyclically symmetric where we begin in order to go around anticlockwise will not matter and this freedom will give rise to identities. A partition $\lambda$ associated with an internal edge will be present in two different triplets and we take the convention that when a partition occurs a second time in a triplet it is replaced with its transpose (this can be understood by considering oriented graphs).

Thus for the above graph we get the factor:
\bea
z^{|\lambda|}\,C_{\emptyset\,\emptyset\,\lambda}(q)\,C_{\emptyset\,\emptyset\,\lambda^{t}}(q)=
z^{|\lambda|}\,C_{\lambda\,\emptyset\,\emptyset}(q)
\,C_{\lambda^{t}\,\emptyset\,\emptyset}(q)=z^{|\lambda|}\,C_{\emptyset\,\lambda\,\emptyset}(q)
\,C_{\lambda^{t}\,\emptyset\,\emptyset}(q)
\eea
where the equality follows from the cyclic symmetry of the vertex.
By summing over all non-trivial partition we get the function associated with the graph:
\bea
{\cal Z}&=&\sum_{\lambda}z^{|\lambda|}\,C_{\emptyset\,\emptyset\,\lambda}(q)\,C_{\emptyset\,\emptyset\,\lambda^{t}}(q)\\\nn
&=&\sum_{\lambda}z^{|\lambda|}\,s_{\lambda^t}(q)\,s_{\lambda}(q)=\prod_{i,j=1}^{\infty}\Big(1+z\,q^{i+j-1}\Big)=\prod_{k=1}^{\infty}\Big(1+z\,q^{k}\Big)^{k}
\eea

Now lets consider the graph obtained from \figref{graph1} by gluing two horizontal external edges as shown in the figure below.
\begin{figure}[h]
\begin{center}
  \includegraphics[width=2.5in]{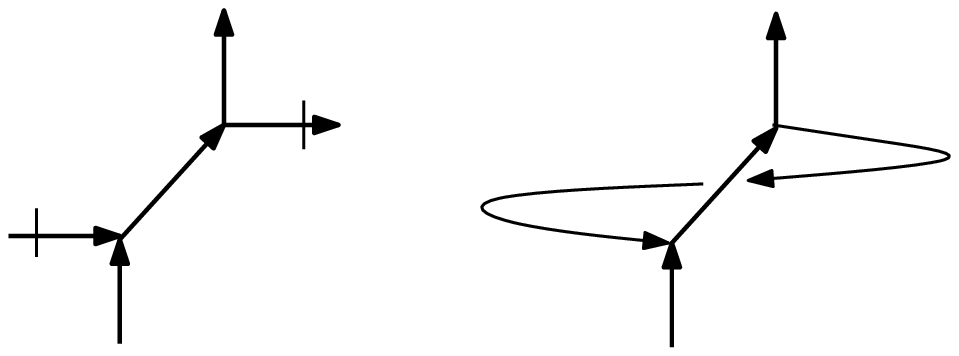}\\
\end{center}
\end{figure}

In this case the function associated with this modified graph is given by
\bea\label{eqw}
{\cal Z}&=&\sum_{\tau,\sigma}z_{1}^{\mid\sigma\mid}z_{2}^{\mid\tau\mid}C_{\emptyset\,
\sigma \,\tau^t}(q)\,C_{\emptyset\,\sigma^t\, \tau}(q)\\\nn
&=&\sum_{\tau ,\sigma }z_{1}^{|\sigma|}\,z_{2}^{|\tau|}\,s_{\tau }(q^{-\rho })\,s_{\tau ^{t}}(q^{-\rho })\,s_{\sigma ^{t} }(q^{-\rho-\tau ^{t}})\,s_{\sigma}(q^{-\rho -\tau })
\eea
Using the summation identity for Schur function,
\bea\label{sum}
\sum_{\mu}q^{\mid\mu\mid}s_{\mu}(x)s_{\mu^{t}}(y)&=&\prod_{i,j\geq1}(1+q\,x_i\,y_j)
\eea
we can sum over $\sigma$ to obtain
\bea\label{hh}
{\cal Z}=\sum_{\tau }z_{2}^{|\tau|}\,s_{\tau }(q^{-\rho })\,s_{\tau ^{t}}(q^{-\rho })\prod_{i,j=1}^{\infty}\Big(1+z_{1}\,q^{-\tau_{i}-\tau^{t}_{j}-\rho_{i}-\rho_{j}}\Big)\,
\eea

The infinite product can be simplified using the identity
\bea
\frac{\prod_{i,j=1}^{\infty}\Big(1+z_{1}\,q^{-\tau_{i}-\tau^{t}_{j}-\rho_{i}-\rho_{j}}\Big)}
{\prod_{i,j=1}^{\infty}\Big(1+z_{1}\,q^{-\rho_{i}-\rho_{j}}\Big)}&=&\prod_{s\in \tau}\Big((1+z_{1}q^{h(s)})(1+z_{1}q^{-h(s)})\Big)\\\nn
&=&z_{1}^{|\tau|}\,q^{-\frac{||\tau||^2+||\tau^{t}||^2}{2}}\,\prod_{s\in \tau}\Big((1+z_{1}q^{h(s)})(1+z_{1}^{-1}q^{h(s)})\Big)
\eea
Using the above identity and the principal specialization of the Schur fucntion,
\bea\label{ps}
s_{\lambda}(q^{-\rho})&=&q^{\frac{\parallel\lambda^t\parallel^{2}}{2}}\prod_{s\in\lambda}(1-q^{h(s)})^{-1}\,,\,\,||\lambda||^2=\sum_{i}\lambda_{i}^2
\eea 
in Eq(\ref{hh}) we get
\begin{equation}
\label{eqww}
{\cal Z}=\Big(\prod_{k\geq 1}(1+z_{1}\,q^{k})^{k}\Big)\sum_{\tau
}(z_{1}z_{2})^{\mid \tau \mid }\prod_{s\in \tau }\frac{%
(1+z_{1}q^{h(s)})(1+z_{1}^{-1}q^{h(s)})}{(1-q^{h(s)})^{2}}
\end{equation}%
Since $C_{\lambda\, \mu\, \nu }$ is cyclically symmetric
\begin{eqnarray}
{C_{\emptyset\, \sigma\, \tau^{t} }(q)} =C_{\sigma\, \tau^{t}\, \emptyset}(q) \,,\,\,\,\,\,\,
{C_{\emptyset \, \sigma ^{t}\, \tau }(q)} ={C_{\sigma ^{t}\, \tau \, \emptyset}(q)}
\end{eqnarray}%
therefore we can write Eq(\ref{eqw}) also as
\bea
{\cal Z}&=&\sum_{\tau\, \sigma }z_{1}^{|\sigma|}\,z_{2}^{|\tau|}\,C_{\sigma\, \tau^{t}\, \emptyset}(q)\,C_{\sigma^{t} \,\tau\,  \emptyset }(q)\\\nn
&=&\sum_{\sigma ,\eta ,\gamma ,\tau }z_{1}^{|\sigma|}\,z_{2}^{|\tau|}\,s_{\sigma ^{t}/\eta }(q^{-\rho })\,s_{\sigma
/\gamma }(q^{-\rho })\,s_{\tau^{t} /\eta }(q^{-\rho })\,s_{\tau /\gamma
}(q^{-\rho })
\eea

Using the summation identities for the skew-Schur function
\bea
s_{\lambda/\mu}(x,y)&=&\sum_{\nu}s_{\lambda/\nu}(x)s_{\nu/\mu}(y)\\\nn
\sum_{\eta}s_{\eta^{t}/\nu}(x)s_{\eta/\mu}(y)&=&\prod_{i,j\geq1}(1+x_iy_j)\sum_{\tau}s_{\mu^{t}/\tau}(x)s_{\nu^{t}/\tau^{t}}(y)\\\nn
\sum_{\rho,\lambda}q^{\mid\rho\mid}s_{\rho/\lambda}(x)s_{\rho/\lambda}(y)&=&\prod_{k\geq1}(1-q^{k})^{-1}\prod_{i,j\geq1}(1-q^k\,x_i\,y_j)^{-1}
\eea

we get
\bea
\label{eqe}
{\cal Z}=\Big(\prod_{k\geq 1}(1+z_{1}\,q^{k})^{k}\Big)\prod_{k\geq
1}(1-z_{1}^{k}z_{2}^{k})^{-1}\prod_{r\geq 1}\frac{
(1+z_{1}^{k-1}z_{2}^{k}q^{r})^{r}(1+z_{1}^{k+1}z_{2}^{k}q^{r})^{r}}{
(1-z_{1}^{k}z_{2}^{k}q^{r})^{2r}}
\eea
Comparing Eq(\ref{eqww}) and Eq(\ref{eqe}) we get the following identity
\begin{eqnarray}\nn
\sum_{\tau }(z_{1}z_{2})^{\mid \tau \mid }\prod_{s\in \tau }\frac{
(1+z_{1}q^{h(s)})(1+z_{1}^{-1}q^{h(s)})}{(1-q^{h(s)})^{2}}=\prod_{k\geq
1}(1-z_{1}^{k}z_{2}^{k})^{-1}\prod_{r\geq 1}\frac{
(1+z_{1}^{k-1}z_{2}^{k}q^{r})^{r}(1+z_{1}^{k+1}z_{2}^{k}q^{r})^{r}}{
(1-z_{1}^{k}z_{2}^{k}q^{r})^{2r}}
\end{eqnarray}%
The above is the generalized form of the Nekrasov-Okounkov identity. If we let
\bea
z_{1}=-e^{\beta\, t}\,,\,\,\,z_{2}=-z\,,\,\, q=e^{-\beta}
\eea
in the above identity then in the limit $\beta\mapsto 0$ we obtain the Nekrasov-Okounkov Identity:
\begin{equation}
{\sum_{\tau }z^{|\tau |}\prod_{s\in \tau }\frac{h^{2}(s)-t^{2}}{%
h^{2}(s)}}={\prod_{k\geq 1}(1-z^{k})^{t^{2}-1}}
\end{equation}

\section{Nekrasov-Okounkov type identities}

To illustrate the method of obtaining these identities we work out another example which leads to another  Nekrasov-Okounkov type identity. The graph we will consider is given in the \figref{graph3} below.

\begin{figure}[h]
\begin{center}
  \includegraphics[width=1.5in]{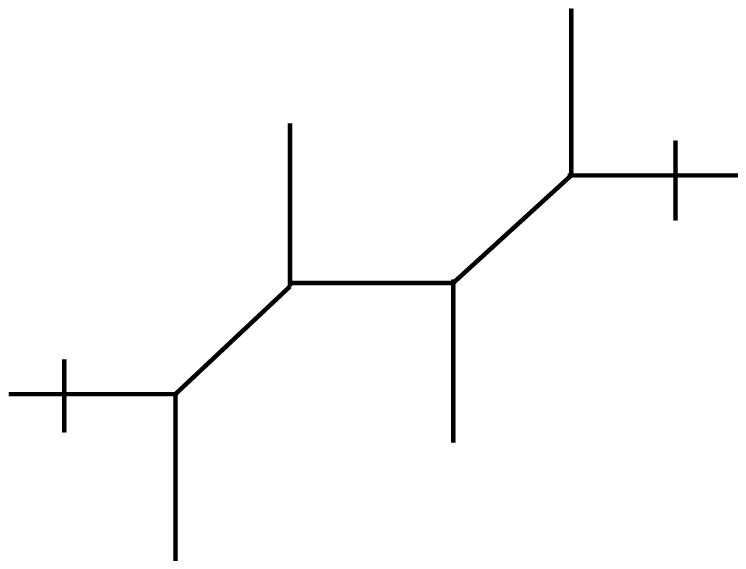}\\
\end{center}
 \label{graph3}
\end{figure}

The function associated with this graph is given by
\bea
\label{ytu}
{\cal Z} &= &\sum _{\tau,\sigma,\nu,\lambda}
z_{1}^{|\sigma|}\,z_{2}^{ |\tau|}\,z_{3}^{|\lambda|}\,z_{4}^ {|\nu|}\,
C_{\sigma^t\,\tau\,\emptyset}(q)\,C_{\sigma\,\nu^t\,\emptyset}(q)\,C_{\lambda^t\,\nu
\,\emptyset}(q)\,C_{\lambda\,\tau^t\,\emptyset}(q)\\\nn
&=&\sum_{\substack{ \tau ,\sigma ,\nu ,\lambda ,  \\ \alpha ,\beta ,\gamma
,\eta }}z_{1}^{|\sigma |}\,z_{2}^{|\tau |}\,z_{3}^{|\lambda
|}\,z_{4}^{|\nu |}s_{{\sigma }/{\eta }}(q^{-\rho })s_{{\tau }/{\eta }%
}(q^{-\rho })s_{{\sigma }^{t}/{\gamma }}(q^{-\rho })  \notag \\\nn
&&s_{{\nu ^{t}}/{\gamma }}(q^{-\rho })s_{{\lambda }/{\beta }}(q^{-\rho })s_{{%
\nu }/{\beta }}(q^{-\rho })s_{{\lambda ^{t}}/{\alpha }}(q^{-\rho })s_{{\tau }%
^{t}/{\alpha }}(q^{-\rho })\\\nn
&=&\prod_{k\geq 1}(1-z^{k})^{-1}\prod_{r\geq 1}\Big(\frac{\prod_{a=1}^{4}(1+z^{k-1}z_{a}q^{r})(1+z^{k}z_{a}^{-1}q^{r})}
{(1-z^{k}q^{r})^{4}(1-\frac{z^{k}}{z_{1}z_{2}}q^{r})
(1-\frac{z^{k}}{z_{2}z_{3}}q^{r})(1-\frac{z^{k}}{z_{3}z_{4}}q^{r})(1-\frac{z^{k}}{z_{4}z_{1}}q^{r})
}\Big)^{r}
\eea
where $z=z_{1}z_{2}z_{3}z_{4}$. Using the cyclic symmetry of the topological vertex we can write the function associated with the graph \figref{graph3} also as
\bea
\label{rer}
{\cal Z}&=&\sum_{\tau ,\sigma ,\nu ,\lambda
}z_{1}^{|\sigma |}\,z_{2}^{|\tau |}\,z_{3}^{|\lambda |}\,z_{4}^{|\nu
|}\,C_{\emptyset\, \sigma ^{t}\,\tau }(q)\,C_{\emptyset\, \sigma\, \nu
^{t}}(q)\,C_{\emptyset\, \lambda^{t}\,\nu}(q)\,C_{\emptyset\, \lambda\, \tau^{t}}(q)\\\nn
&=&\sum_{\tau ,\sigma ,\nu ,\lambda
}z_{1}^{|\sigma |}\,z_{2}^{|\tau |}\,z_{3}^{|\lambda |}\,z_{4}^{|\nu
|}\,s_{\sigma }(q^{-\rho -\nu })s_{\sigma ^{t}}(q^{-\rho -\tau
^{t}})s_{\lambda ^{t}}(q^{-\rho -\nu ^{t}})s_{\lambda }(q^{-\rho -\tau
})s_{\tau }(q^{-\rho })s_{\tau ^{t}}(q^{-\rho })s_{\nu }(q^{-\rho })s_{\nu
^{t}}(q^{-\rho })  \\\nn
\eea
Summation over $\sigma$ and $\lambda$ can be carried out easily using Eq(\ref{sum}) and using the infinite product identity 
\bea\nn
\frac{\prod_{i,j=1}^{\infty}\Big(1+z_{1}q^{-\nu_{i}-\tau^{t}_{j}+i+j-1}\Big)}{\prod_{i,j=1}^{\infty}\Big(1+z_{1}q^{i+j-1}\Big)}=
\prod_{(i,j)\in \nu}\Big(1+z_{1}\,q^{-\nu_{j}-\tau_{i}^{t}+i+j-1}\Big)\,\prod_{(i,j)\in \tau}\Big(1+z_{1}\,q^{\nu_{j}^{t}+\tau_{i}-i-j+1}\Big)
\eea
we get
\bea\label{id}
{\cal Z} &=&\Big(\prod_{r\geq
1}(1+z_{1}q^{r})^{r}(1+z_{3}q^{r})^{r}\Big)\sum_{\tau \,,\nu }(z_{1}z_{4})^{|\tau
|}(z_{2}z_{3})^{|\nu |} \\\nn
&&\prod_{(i,j)\in \tau }\frac{(1+z_{1}^{-1}q^{i+j-1-\nu _{j}^{t}-\tau
_{i}})(1+z_{3}q^{i+j-1-\nu _{j}^{t}-\tau _{i}})}{(1-q^{h(i,j)})^{2}}
\,\prod_{(i,j)\in \nu }\frac{(1+z_{1}q^{i+j-1-\nu _{j}-\tau
_{i}^{t}})(1+z_{3}^{-1}q^{i+j-1-\nu _{j}-\tau _{i}^{t}})}{%
(1-q^{h(i,j)})^{2}}
\eea
Comparing Eq(\ref{ytu}) and Eq(\ref{id}) we get the identity
\begin{eqnarray}
&&\sum_{\tau \,,\nu }(z_{1}z_{4})^{|\tau
|}(z_{2}z_{3})^{|\nu |} \,\prod_{(i,j)\in \tau }\frac{(1+z_{1}^{-1}q^{i+j-1-\nu _{j}^{t}-\tau
_{i}})(1+z_{3}q^{i+j-1-\nu _{j}^{t}-\tau _{i}})}{(1-q^{h(i,j)})^{2}}\\\nn
&&\,\prod_{(i,j)\in \nu }\frac{(1+z_{1}q^{i+j-1-\nu _{j}-\tau
_{i}^{t}})(1+z_{3}^{-1}q^{i+j-1-\nu _{j}-\tau _{i}^{t}})}{%
(1-q^{h(i,j)})^{2}}\\\nn
&=&\prod_{k\geq 1}(1-z^{k})^{-1}\prod_{r\geq 1}\Big(\frac{(1+z_{2}q^{r})(1+z_{4}q^{r})\prod_{a=1}^{4}(1+z^{k}z_{a}q^{r})(1+z^{k}z_{a}^{-1}q^{r})}
{(1-z^{k}q^{r})^{4}(1-\frac{z^{k}}{z_{1}z_{2}}q^{r})
(1-\frac{z^{k}}{z_{2}z_{3}}q^{r})(1-\frac{z^{k}}{z_{3}z_{4}}q^{r})(1-\frac{z^{k}}{z_{4}z_{1}}q^{r})
}\Big)^{r}
\end{eqnarray}
Let us define
\bea\label{sss}
z_{1}=-e^{\beta\, t_{1}}\,,\,\,\,z_{3}=-e^{\beta\,t_{3}}\,,z_{2}=-x\,,\,\,\,z_{4}=-y\,,\,\,\,\,\,\,q=e^{-\beta}\,
\eea
then in the limit $\beta \rightarrow 0$ the above identity reduces to
\begin{eqnarray}
&&\sum_{\tau ,\nu }x^{|\nu|}\,y^{|\tau|}\prod_{s\in \tau }%
\frac{(\ell_{\tau}(s)+a_{\nu}(s)+1+t_{3})(\ell_{\tau}(s)+a_{\nu}(s)+1-t_{1}))}{h^{2}(s)}  \\\nn
&&\prod_{s\in \nu^{t} }\frac{(\ell_{\nu}(s)+a_{\tau}(s)+1+t_{3})(\ell_{\nu}(s)+a_{\tau}(s)+1-t_{1}))}{h^{2}(s)}  \\\nn
&=&\prod_{k\geq
1}(1-x^{k}\,y^{k-1})^{t_{1}t_{3}}(1-x^{k-1}\,y^{k})^{t_{1}t_{3}}(1-x^{k}\,y^{k})^{t_{1}^{2}+t_{3}^{2}-1}
\end{eqnarray}

\section{Conclusions}
We have shown that Nekrasov-Okounkov identity and its generalizations can be obtained using the cyclic symmetry of the topological vertex. This method can be used to generate highly non-trivial identities involving Schur and Skew-Schur functions. These identities encode in them certain dualities of the underlying physical theories. It would be interesting to prove these identities using purely combinatorial methods as was carried out for the Nekrasov-Okounkov identity in \cite{han}. This method also sheds light on a question of Richard Stanley. In \cite{stanley} Stanley asked if the Nekrasov-Okounkov identity Eq(\ref{NO}) can be generalized such that left hand side is
\bea\label{ss}
{\sum_{\lambda}z^{|\lambda|}\,\,\prod_{s\in \lambda}\prod_{a=1}^{k}\frac{
h(s)^2-t_{a}^2}{h(s)^2}}
\eea
Using the method discussed before we can obtain the following identity:
\bea
&&\sum_{\nu}(z\,x\,y)^{|\nu|}\prod_{s\in\nu}\frac{(1+x\,q^{h(s)})(1+x^{-1}q^{h(s)})(1+y\,q^{h(s)})(1+y^{-1}q^{h(s)})}{(1-q^{h(s)})^4}
\\\nn
&=&\sum_{\nu,\mu,\lambda,\eta_{1},\eta_{2},\eta_{3}\eta_{4}}
z^{|\nu|}x^{|\mu|}y^{|\lambda|}\,
s_{\nu/\eta_{1}}({\bf x})s_{\nu^{t}/\eta_{2}}({\bf x})
s_{\nu/\eta_{3}}({\bf x})s_{\nu^{t}/\eta_{4}}({\bf x})
s_{\mu/\eta_{2}}({\bf x})s_{\mu^{t}/\eta_{3}}({\bf x})
s_{\lambda/\eta_{4}}({\bf x})s_{\lambda^{t}/\eta_{1}}({\bf x})
\eea
where ${\bf x}=q^{-\rho}=\{q^{1/2},q^{3/2},q^{5/2},\cdots\}$. The left hand side of this identity can indeed be reduced to Eq(\ref{ss}) (for $k=2$) by taking a limit similar to the one in Eq(\ref{sss}) but the right hand side involves products of skew-Schur functions. We are not familiar with any identities which will allow such a product of skew-Schur function to be summed up in to a product. It would be interesting to find a product representation of this sum. The case $k>2$ can be worked out as well and gives again similar sum over the product of skew-Schur functions.

\end{document}